\documentclass[a4paper, 11pt]{article}

\usepackage{geometry}

\usepackage{amsfonts}
\usepackage{amssymb}
\usepackage{amsmath}
\usepackage{amsthm}
\usepackage{verbatim}
\usepackage{enumitem}
\usepackage{xcolor}

\theoremstyle{plain}
\newtheorem{theorem}{Theorem}[section]

\newtheorem{observation}[theorem]{Observation}
\theoremstyle{definition}

\usepackage{doc}

\title{Game-theoretic semantics and 
partial specifications
}

\author{
Antti Kuusisto\\
Tampere University and University of Helsinki
}

\date{}

\begin{document}

\maketitle

\begin{abstract}
\noindent
\noindent
We discuss partial specifications in
first-order logic FO and also in a Turing-complete extension of FO.
We compare the compositional and game-theoretic approaches to
the systems.
\end{abstract}


\setcounter{tocdepth}{2}

%
%


\section{Introduction}
A natural Turing-complete extension $\mathcal{L}$ of
first-order logic $\mathrm{FO}$ can be obtained by extending $\mathrm{FO}$ by
two new features: (1)
operators that \emph{add new points} to models
and \emph{new tuples} to relations, and
(2) operators that enable formulae to \emph{refer to themselves}.
The self-referentiality operator of $\mathcal{L}$ is based on a construct that enables \emph{looping}
when formulae are evaluated using game-theoretic semantics.
The reason the logic $\mathcal{L}$ is
particularly interesting lies in its simplicity and its
\emph{exact behavioural
correspondence} with Turing machines.
It is also worth noting that the new operators of $\mathcal{L}$ nicely
capture two fundamental classes of constructors that are omnipresent in
the everyday practice of mathematics (but missing from $\mathrm{FO}$): 
scenarios where fresh points are added to investigated constructions
(or fresh lines are drawn, et cetera) play a central role in geometry, and
recursive looping operators are found everywhere in mathematical practice,
often indicated with the help of the famous three dots (...).

Adding new points to domains makes functions partial, and thus it is 
natural to consider partial functions in the framework of $\mathcal{L}$.
Previous work has concentrated on relational vocabularies. We develop a
related semantics for $\mathcal{L}$, also accommodating partial relation symbols
into the picture. We begin by investigating a related system of first-order logic 
via a compositional semantics and then extend it to a game-theoretic semantics for $\mathcal{L}$.

The logic $\mathcal{L}$ was first introduced in \cite{kuusistoturingcomp}.
Below we discuss the results of that article and also further results not covered there.
Other systems that bear some degree of similarity to $\mathcal{L}$ include for example BGS logic 
\cite{bgs} and abstract state machines
\cite{guuure1, borger}.
%
%
Logics that modify structures include, for example, sabotage modal logic and
public announcements logics. Recursive looping operators are a common feature in
logics in finite model theory and verification. However, while the approach in $\mathcal{L}$ bears a
degree of similarity to the fixed point operators of the $\mu$-calculus, $\mathcal{L}$ is not based on
fixed points and no monotonicity requirements apply.

\section{Partially fixed entities}

Let $\tau$ be a vocabulary, i.e., it is a collection of relation symbols, function 
symbols and constant symbols. Let $\mathrm{VAR}$ denote a
countably infinite set of first-order variable symbols. Define $\tau$-terms in the usual way to be
the smallest set of terms built from variables $x\in \mathrm{VAR}$ and constant symbols $c\in \tau$ by
applying the function symbols $f\in \tau$ in the way that respects the function symbol arities, and
let $\mathrm{FO}(\tau)$ denote the language of first-order logic when the vocabulary under 
investigation is $\tau$. Now, let $\mathfrak{M}$ be a $\tau$-model where the function 
symbols $f\in \tau$ can be interpreted as partial functions on the model domain $M$. A $k$-ary function
symbol $f\in \tau$ is interpreted so that $f^{\mathfrak{M}}$ is a
function $f^{\mathfrak{M}}:N\rightarrow M$ where $N$ is some subset of $M^k$.
We allow also the constant symbols to be interpreted as partial, meaning that $c\in \tau$ can be interpreted as an
element of $\mathfrak{M}$ or it may not have an interpretation in $\mathfrak{M}$ at all.

In addition to the model $\mathfrak{M}$, we also need a \emph{variable assignment} $g$,
which is a function $g:V\rightarrow M$ mapping from some set $V\subseteq\mathrm{VAR}$
into the domain $M$ of $\mathfrak{M}$.
Now consider a $\tau$-term $t$. If $t = f$,
then $t^{\mathfrak{M},g}$ is the partial function $f^{\mathfrak{M},g}$.  We may 
also write $f^{\mathfrak{M}}$.
When $t= x\in \mathrm{VAR}$, then $t^{\mathfrak{M},g} = x^{\mathfrak{M},g} = g(x)$ when $x$ belongs to
the domain of the variable assignment $g$, and
otherwise we say that $t^{\mathfrak{M},g} = x^{\mathfrak{M},g}$
does not exist (the variable is undefined). When $t = c \in \tau$ is a constant symbol, 
then $t^{\mathfrak{M},g} = c^{\mathfrak{M},g}$ is 
the element $u\in M$ if $\mathfrak{M}$ specifies such an element $u$ for $c$, and
otherwise we say that $t^{\mathfrak{M},g} = c^{\mathfrak{M},g}$ is undefined. 
If $t = f(t_1,\dots , t_n)$, then $t^{\mathfrak{M},g}
= f^{\mathfrak{M}}(t_1^{\mathfrak{M},g},\dots , t_n^{\mathfrak{M},g})$ if 
each of $t_1^{\mathfrak{M},g},\dots , t_n^{\mathfrak{M},g}$ is defined and $f^{\mathfrak{M}}$ is
defined on the tuple $(t_1^{\mathfrak{M},g},\dots , t_n^{\mathfrak{M},g})$; otherwise we
say that $t^{\mathfrak{M},g}$ is not defined.
Recall that $g[x\mapsto a]$ is the assignment 
with domain $\mathit{dom}(g)\cup \{x\}$ that is otherwise as $g$
but maps $x$ to $a$.
The semantics of first-order logic $\mathrm{FO}(\tau)$ is 
defined as follows.

\[
\begin{array}{lll}
%
%
%
%
\mathfrak{M}, g\models^+ R(t_1,\dots , t_n) &
{\Leftrightarrow} & \text{ the terms }
          t_1^{\mathfrak{M},g},\dots , t_n^{\mathfrak{M},g}
\text{ are defined and }\\
& &\text{ and we have }
(t_1^{\mathfrak{M},g},\dots , t_n^{\mathfrak{M},g})\in R^{\mathfrak{M}}\\
\mathfrak{M}, g\models^+ t_1 = t_2 &
{\Leftrightarrow} & t_1^{\mathfrak{M},g}\text{ and } t_2^{\mathfrak{M},g}
\text{ are defined and }t_1^{\mathfrak{M},g} =  t_2^{\mathfrak{M},g}\\
\mathfrak{M}, g\models^+ \neg\varphi &
{\Leftrightarrow} & \mathfrak{M}, g \models^- \varphi \\
\mathfrak{M}, g\models^+ \varphi\wedge \psi & {\Leftrightarrow} &
\mathfrak{M}, g\models^+ \varphi \text{ and }\mathfrak{M}, g\models^+ \varphi \\
\mathfrak{M}, g\models^+ \exists x\, \varphi& {\Leftrightarrow} &
\text{there exists }a\in M\text{ such that }
\mathfrak{M}, g[x \mapsto a]\models^+  \varphi\\
%
%
%
%
\mathfrak{M}, g\models^- R(t_1,\dots , t_n) &
{\Leftrightarrow} & \text{the terms }
          t_1^{\mathfrak{M},g},\dots , t_n^{\mathfrak{M},g}
\text{ are defined and }\\
& &\text{and we have }
(t_1^{\mathfrak{M},g},\dots , t_n^{\mathfrak{M},g})\not\in R^{\mathfrak{M}}\\
\mathfrak{M}, g\models^- t_1 = t_2 &
{\Leftrightarrow} & t_1^{\mathfrak{M},g}\text{ and } t_2^{\mathfrak{M},g}
\text{ are defined and }t_1^{\mathfrak{M},g} \not=  t_2^{\mathfrak{M},g}\\
\mathfrak{M}, g\models^- \neg\varphi &
{\Leftrightarrow} & \mathfrak{M}, g \models^+ \varphi \\
\mathfrak{M}, g\models^- \varphi\wedge \psi & {\Leftrightarrow} &
\mathfrak{M}, g\models^- \varphi \text{ or }\mathfrak{M}, g\models^- \varphi \\
\mathfrak{M}, g\models^- \exists x\, \varphi& {\Leftrightarrow} &
\text{for all }a\in M\text{ we have }
\mathfrak{M}, g[x \mapsto a]\models^-  \varphi\\
\end{array}
\]

\medskip

Note that for a constant symbol $c$, we can write $c^{\mathfrak{M}}$ as
well as $c^{\mathfrak{M},g}$.

We can also consider partially defined relations. A partially defined
relation $R^{\mathfrak{M}}$ specifies (1) a set of tuples $(u_1,\dots , u_n)$ in the
relation and (2) a set of tuples $(v_1,\dots , v_n)$ \emph{not} in the
relation. Thus some tuples may not be specified either way (neither positive 
nor negative) and thus it is 
undefined whether such a tuple belongs to the relation.
In such a framework, we get the modified clauses

\[
\begin{array}{lll}
%
%
%
%
\mathfrak{M}, g\models^+ R(t_1,\dots , t_n) &
{\Leftrightarrow} & \text{the terms }
\text{ are defined and }\\
& & R^{\mathfrak{M}}\text{ is
defined on }(t_1^{\mathfrak{M},g},\dots , t_n^{\mathfrak{M},g})\\
& &\text{and we have }
(t_1^{\mathfrak{M},g},\dots , t_n^{\mathfrak{M},g})\in R^{\mathfrak{M}}\\
%

%
%
%
%
\mathfrak{M}, g\models^- R(t_1,\dots , t_n) &
{\Leftrightarrow} & \text{the terms }
          t_1^{\mathfrak{M},g},\dots , t_n^{\mathfrak{M},g}
\text{ are defined and }\\
& & R^{\mathfrak{M}}\text{ is
defined on }(t_1^{\mathfrak{M},g},\dots , t_n^{\mathfrak{M},g})\\
& &\text{and we have }
(t_1^{\mathfrak{M},g},\dots , t_n^{\mathfrak{M},g})\not\in R^{\mathfrak{M}}\\
\end{array}
\]

\medskip

It is also possible to add a contradictory negation $\sim$ into the picture.
We can let, to give a possibility,

\medskip

\[
\begin{array}{lll}
\mathfrak{M}, g\models^+ \sim\varphi &
{\Leftrightarrow} & \text{not }\mathfrak{M}, g \models^+ \varphi\\
\mathfrak{M}, g\models^- \sim\varphi &
{\Leftrightarrow} & \text{not }\mathfrak{M}, g \models^- \varphi.\\
\end{array}
\]

We also define the determinacy operator $d$ such that

\[
\begin{array}{lll}
\mathfrak{M}, g\models^+ d\, \varphi &
{\Leftrightarrow} & \mathfrak{M}, g \models^+ \varphi
\text{ or }\mathfrak{M}, g \models^- \varphi\\
\mathfrak{M}, g\models^- d\, \varphi &
{\Leftrightarrow} & (\text{not }\mathfrak{M}, g \models^+ \varphi)
\text{ and }(\text{not }\mathfrak{M}, g \models^- \varphi)\\
\end{array}
\]

We note that further generalizations of indeterminate 
values come naturally. For example, one could consider $k$-ary
function symbols that are defined on some $(k+1)$-tuple $(u_1,\dots , u_{k+1})$
such that $f^{\mathfrak{M}}(u,\dots , u_k)$ is
defined to \emph{not} equal $u_{k+1}$.
However, we shall not consider such extensions here.

\section{The logic $\mathcal{L}$}

We then consider the Turing complete logic $\mathcal{L}$, or \emph{computation logic}
$\mathrm{CL}$,
sometimes also called \emph{computation game logic}.
Let us first consider a relational vocabulary with no constant or function symbols.
Let $\mathcal{L}$ denote the language that extends the
syntax of first-order logic by the following formula construction rules:
%

\vspace{-2mm}

\begin{enumerate}[noitemsep]
\item
$\varphi\ \mapsto\ \mathrm{I}x\, \varphi$
\item
$\varphi\ \mapsto\ \mathrm{I}_{R(x_1, \dots , x_n)}\ \varphi$
\item
$\varphi\ \mapsto\ \mathrm{D}x\, \varphi$
\item
$\varphi\ \mapsto\ \mathrm{D}_{R(x_1, \dots , x_n)}\ \varphi$
\item
$C_i$ is an atomic formula (for each $i\in \mathbb{N}$)
\item
$\varphi\ \mapsto\ \, C_i\, \varphi$
\item
We also allow atoms $X(x_1,\dots , x_n)$ where $X$ is a
relational symbol not in the vocabulary considered. These symbols are
analogous to tape symbols. The simplest way to treat $X$ is to consider it a relation
symbol interpreted initially as the empty $n$-ary relation $\emptyset$.
\end{enumerate}
Intuitively, a formula of type $\mathrm{I}x\, \varphi(x)$ states that it is 
possible to \emph{insert} a fresh, isolated element $u$ into the domain of the current model so
that the resulting new model satisfies $\varphi(u)$. The fresh element $u$ being \emph{isolated} 
means that $u$ is disconnected from the original model; the relations of the original model are
not altered in any way by the operator $\mathrm{I}x$, so $u$ does not become part of
any relational tuple at the moment of insertion.

A formula of type $\mathrm{I}_{R(x_1, \dots , x_n)}\ \varphi(x_1,\dots , x_n)$
states that it is possible to insert a tuple $(u_1,\dots , u_n)$ into the relation $R$ so
that $\varphi(u_1,\dots , u_n)$ holds in the obtained model. The tuple $(u_1,\dots , u_n)$ is a
sequence of elements in the original model, so this time the domain of the model is 
not altered. Instead, the $n$-ary relation $R$ obtains a new
tuple via the insertion. 
A formula of type $\mathrm{D}x\, \varphi$ states that it is 
possible to \emph{delete} an element $u$ named $x$ from the domain of the current model so
that the resulting new model satisfies $\varphi$. All tuples that contain $u$ are of course also
deleted from the related relations.
A formula of type $\mathrm{D}_{R(x_1, \dots , x_n)}\ \varphi$
states that we can delete a tuple $(u_1,\dots , u_n)$ 
named $(x_1,\dots , x_n)$ from the relation $R$ so
that $\varphi$ holds in the obtained model. 
The new atomic formulae $C_i$ can be regarded as \emph{variables} ranging over
formulae, so a formula $C_i$ can be considered to be a \emph{pointer} to (or the
\emph{name} of) some other formula.
The formulae $C_i\, \varphi$ could intuitively be given the following reading:
\emph{the claim $C_i$, which states that $\varphi$, holds.}
Thus the formula $C_i\, \varphi$ is both \emph{naming} $\varphi$ to be called $C_i$
and \emph{claming} that $\varphi$ holds.
Importantly, the formula $\varphi$ can contain $C_i$ as an atomic formula.
This leads to self-reference.

The logic $\mathcal{L}$ is based on game-theoretic semantics GTS which 
directly extends the standard GTS of $\mathrm{FO}$. Recall that the GTS of $\mathrm{FO}$ is
based on games played by \emph{Eloise} (who is initially the \emph{verifier}) and 
\emph{Abelard} (initially the \emph{falsifier}). In a 
game $G(\mathfrak{M},g,\varphi)$, the verier is trying to show (or 
verify) that $\mathfrak{M},g\models\varphi$ and the falsifier is opposing this, i.e.,
the falsifier wishes to falsify the claim $\mathfrak{M},g\models\varphi$. The 
players start from the \emph{position} $(\mathfrak{M},g,\varphi)$ and work their way towards the
subformulae of $\varphi$. See $\cite{mann}$ for a detailed exposition of GTS for
first-order logic and \cite{hintikka22} for some of the founding ideas behind GTS.

To deal with the logic $\mathcal{L}$,
the rules for the $\mathrm{FO}$-game are extended as follows. 

\vspace{-1mm}

\begin{enumerate}[noitemsep]
\item
In a position $(\mathfrak{M},g,\mathrm{I}x\, \psi)$, the game is
continued from a position $(\mathfrak{M}',g[x\mapsto u],\psi)$ where $\mathfrak{M}'$ is
the model obtained by inserting a fresh isolated point $u$ into the 
domain of $\mathfrak{M}$.
\item
In a position $(\mathfrak{M},g,\mathrm{I}_{R(x_1,\dots , x_n)}\, \psi)$, 
the verifier chooses a tuple $(u_1,\dots , u_n)$ of elements in $\mathfrak{M}$ and
the game is continued from the position
$$(\mathfrak{M}',g[x_1\mapsto u_1,
\dots , x_n\mapsto  u_n],\psi)$$
where $\mathfrak{M}'$ is
the model obtained from $\mathfrak{M}$ by inserting the tuple $(u_1,\dots , u_n)$ into
the relation $R$. If the model has an empty domain (and thus no tuple can be inserted to a
relation of positive arity),
the game play ends and the verfier loses the play.
\item
In a position $(\mathfrak{M},g,\mathrm{D}x\, \psi)$, the game is
continued from a position $(\mathfrak{M}',g^*,\psi)$ where $\mathfrak{M}'$ is
the model obtained by deleting the point $g(x)$ (if it exists) from the
domain of $\mathfrak{M}$. The assignment $g^*$ is obtained 
from $g$ by deleting the pairs $(y,u)$ where $u = g(x)$.
If  $x$ does not belong to the domain of $g$, then the game ends and
the verifier loses the play of the game.\footnote{An alternative
convention would be to just ignore the deletion when it is not possible.}
\item
In a position $(\mathfrak{M},\mathrm{D}_{R(x_1,\dots , x_n)}\, \psi)$, 
the verifier chooses a tuple $(u_1,\dots , u_n)$ of elements in $\mathfrak{M}$ and
the game is continued from the position
$$(\mathfrak{M}',g[x_1\mapsto u_1,
\dots , x_n\mapsto u_n],\psi)$$
where $\mathfrak{M}'$ is
the model obtained from $\mathfrak{M}$ by deleting the tuple $(g(x_1),\dots , g(x_n))$ from
the relation $R$. If some $x_i$ does not belong to the domain of $g$, then the game ends and
the verifier loses the play of the game.\footnote{An alternative
convention would again be to just ignore the deletion when it is not possible.}
\item
In a position  $(\mathfrak{M},g,C_i\, \psi)$, the game
simply moves to the position $(\mathfrak{M},g,\psi)$.
\item
In an atomic position  $(\mathfrak{M},g,C_i)$, the game
moves to the position $(\mathfrak{M},g,C_i\, \psi)$. Here $C_i\, \psi$ is a subformula of
the original formula $\varphi$ that the semantic game began with. (If there are many such
subformulae, the verifier chooses which one to continue from. If there are none, 
the game ends and neither player wins that play of the game. Alternative conventions
would be possible here as well, like the verifier losing.)
\item
In a position with $\wedge,\neg,\exists$, the game proceeds as in first-order logic;
see also \cite{kuusistoturingcomp} which discusses the logic $\mathcal{L}$.
We of course need an extra flag in positions to account for the negation $\neg$, et cetera.
\end{enumerate}

Just like the semantic game for first-order logic, the extended game ends if an atomic
position with a first-order atom $R(x_1,\dots , x_n)$ or $x=y$ is encountered.
The winner is then decided precisely as in the $\mathrm{FO}$-game.\footnote{If some
variable of an atom is not interpreted by $g$, then neither player wins.} Thus the 
extended game can go on forever, as for example the games for $C_i\, C_i$ and $C_i\, \neg\, C_i$
always will. In the case the play of the game indeed goes
on forever, then that play is won by \emph{neither} of the players. Note that Turing-machines 
exhibit precisely this kind of behaviour: they can

\begin{enumerate}[noitemsep]
\item
\emph{halt in an accepting state} 
(corresponding to the verifier winning the semantic game play),
\item
\emph{halt in a 
rejecting state} (corresponding to the falsifier winning),
\item
\emph{diverge} (corresponding to neither of the players winning).
\end{enumerate}

Indeed, there  is an exact correspondence
between the logic $\mathcal{L}$ and Turing machines.
Let $\mathfrak{M},g\models^+\varphi$ (respectively, $\mathfrak{M},g\models^-\varphi$)
denote that Eloise (respectively, Abelard) 
has a winning strategy in the game starting from $(\mathfrak{M},g,\varphi)$.
We may drop $g$
when it is the empty assignment.
Let $\mathrm{enc}(\mathfrak{M})$
denote the encoding of the \emph{finite} model $\mathfrak{M}$ according to any
standard encoding scheme.
Then the following theorem shows that $\mathcal{L}$ corresponds to Turing machines so that not
only acceptance and rejection but even divergence of Turing computation is captured in a
precise and natural way. The proof of the following theorem 
follows from \cite{kuusistoturingcomp}.

\noindent
\begin{theorem}\label{theorem1}
For every Turing machine $\mathrm{TM}$, there
exists a formula $\varphi\in\mathcal{L}$
such that
\begin{enumerate}[noitemsep]
\item
$\mathrm{TM}$ {accepts} $\mathrm{enc}(\mathfrak{M})$\hspace{1mm}
iff\hspace{3mm} $\mathfrak{M}\, {\models}^+\, \varphi$
\item
$\mathrm{TM}$ {rejects} $\mathrm{enc}(\mathfrak{M})$\hspace{1mm}
iff\hspace{1mm} $\mathfrak{M}\, {\models^-}\, \varphi$
\end{enumerate}
Vice versa, for every formula $\varphi\in\mathcal{L}$, there is a
Turing machine $\mathrm{TM}$ such that the above conditions hold.
\end{theorem}

We briefly digress to discuss the level of naturalness of the features of $\mathcal{L}$.
Since $\mathcal{L}$ defines precisely 
the recursively enumerable classes of finite models, it cannot be closed 
under negation (meaning complement here).
Thus $\neg$ is not the classical negation. However, $\mathcal{L}$
has a very natural translation into natural language. The key is to
replace \emph{truth} by \emph{verification}. We read $\mathfrak{M}\models^+\varphi$ as
the claim that ``\emph{it is verifiable that $T(\varphi)$}" where $T$ is
the translation from $\mathcal{L}$ into natural language defined as follows.
We let $T$ map $\mathrm{FO}$-atoms in the
usual way to the corresponding natural language statements, so for
example $T(x=y)$ simply reads ``$x$ equals $y$".
The atoms $C_i$ are read as they stand, so $T(C_i) = C_i$.
The $\mathrm{FO}$-quantifiers translate in the standard way, so 
$T(\exists x \varphi) = \text{\emph{there exists an} $x$ 
\emph{such that} }T(\varphi)$ and analogously for $\forall x$.
Also $\vee$ and $\wedge$
translate in the standard way, so
$T(\varphi\vee \psi) = T(\varphi)\text{ \emph{or} }T(\psi)$ and analogously for $\wedge$.
However, $T(\neg\psi) = \text{\emph{it is falsifiable that }}T(\psi)$ (or alternatively, $\text{\emph{it is refutable that }}T(\psi)$).  Thus negation now translates to
the dual of verifiability. Concerning the
insertion operators, we let $T(\mathrm{I}x\, \varphi) = \emph{it is possible to insert a 
new element}$ $ \emph{$x$ such}$ $\emph{that }T(\varphi)$. Similarly, we let 
$T(\mathrm{I}_{R(x_1,\dots , x_n)}\, \varphi)
= \emph{it is possible to insert a tuple } (x_1,\dots , x_n)$
\emph{into $R$ such that }$T(\varphi)$. Finally, we 
let \[T(C_i\, \varphi\, )
= \emph{it is possible verify the claim } C_i \emph{ which states that }  T(\varphi).\]
Alternatively, we can let \[T(C_i\, \varphi\, )
= \emph{it is possible verify the claim named } C_i \emph{ which claims that }  T(\varphi).\]
Deletion operators are similar to the insertion operators.

Thereby the logic $\mathcal{L}$ can be seen as a \emph{simple Turing-complete fragment of
natural language}. Thus it is not just a technical logical formalism.

It is interesting to note that $\neg$ can---and perhaps should---be 
read as the classical negation on those fragments of $\mathcal{L}$ where the 
semantic games are determined (such as standard $\mathrm{FO}$ contexts).
Furthermore, adding a generalized quantifier to $\mathcal{L}$ 
corresponds precisely to adding a corresponding oracle to Turing machines.

Theorem \ref{theorem1} holds also without tape predicates, given 
the underlying vocabulary contains at least one binary (or
higher-arity) relation.

\begin{observation}
The claim of Theorem \ref{theorem1} holds for $\mathcal{L}$ without 
tape predicates in the scenario where the underlying vocabulary 
contains at least one relation symbol of arity at least two. 
\end{observation}

\begin{proof}
The argument is based on using gadgets. 
Consider a formula $\varphi$ that makes use of 
tape predicates. To eliminate these, we 
will write a new formula $\varphi'$. We assume, 
without loss of generality, that the underlying 
vocabulary contains a binary relation $R$. In 
the case there is no binary relation, we can easily
modify our argument to account for that essentially by 
using the first two coordinate positions of 
some higher-arity relation to
encode a binary relation.

We begin $\varphi'$
with $Iz_1$ which introduces a new domain element $m_1$ labeled by $z_1$. 
Similarly, we use 
further insertion operators to
construct a fresh \emph{successor structure} with 
the elements $m_1, \dots , m_k$, 
where $k$ is the number of tape predicates used in $\varphi$. 
By a successor structure, we mean that $R$ acts as a successor 
relation over the new elements, thereby
connecting the elements $m_1,\dots , m_k$ in the 
given order. 
This way we obtain a fresh point $m_i$ for each
tape predicate $i$ in $\varphi$, and we
have a successor ordering (via $R$) of the fresh points. Each point $m_i$ is
labeled with the corresponding variable $z_i$.

The next task is show how to construct (or model) tuples of the tape predicate $i$. 
This is done via the following reification technique. 
Let $A$ be the domain of the current model minus all the 
fresh points (such as $m_1,\dots , m_k$) that we
will use for encoding tuples of tape predicates. We
call $A$ the proper domain of the model. An
arbitrary $n$-ary tuple $(a_1,\dots , a_n) \in A^n$ of the predicate $i$
will be encoded using a new fresh domain point $u\not\in A$ (called \emph{tuple point})
such that $(m_i,u) \in R$ and such that the 
following conditions hold.

\begin{enumerate}
\item
For each $p\in \{1,\dots , n\}$, 
there exist a fresh point $v_p\not\in A$ (called \emph{coordinate point})
outside the 
proper domain $A$  
such that $(u,v_p)\in R$. 
\item
$(v_p,a_p) \in R$ holds for each $p\in \{1,\dots , n\}$. 
\item
There is a successor order using $R$ and 
ordering the points $v_p$, i.e., $(v_i,v_{i+1})\in R$
for each $i\in \{1, \dots , n-1\}$. This is
essential for distinguishing the order of
the elements of the tuple $(a_1,\dots , a_n)$ 
being encoded. 
\end{enumerate}
In summarly, if the points $m_i$ are 
called \emph{predicate points}, then a predicate 
point connects to a tuple point, a tuple 
point to a coordinate point, and
finally a coordinate point to
the a point of the tuple (over
the proper domain) being encoded. 
Furthermore, we order the predicate 
points as well as the 
coordinate points by a successor order.

It is easy to see that this kind of an encoding can be 
used to eliminate tape predicates. The underlying 
language is strong enough in expressive 
power to deal with also the updates of the 
predicate encodings during semantic games. 

Finally, if nullary tape predicates are used, the way to
encode them is to first turn them into unary tape 
predicates with an appropriate translation and
then encode the obtained unary predicates in the above way. 
\end{proof}

We note that the above also proves the corresponding statement 
for the logic $\mathcal{L}$ as it was defined in
\cite{kuusistoturingcomp}. That is, let $\mathcal{L}_t$ be
the tape-predicate-free fragment of the logic $\mathcal{L}$ as
defined \cite{kuusistoturingcomp}.
Then Theorem \ref{theorem1} 
and similarly Theorems 5.1 and 5.2 of \cite{kuusistoturingcomp} hold
for $\mathcal{L}_t$. The same applies to other close variants of 
the logic, e.g., ones that treat the semantics of variable symbols
without interpretation in a slightly modified way. The above observations
were essentially made already in \cite{kuusistoturingcomp}, but left
there officially as conjectures for the sake of lack of space.

\subsection{Partially fixed entities in $\mathcal{L}$}

The logic $\mathcal{L}$ has, as such, natural indeterminacy properties
stemming from divergence of Turing machines and also in other ways.
The truth teller formula $C_i\, C_i$ is indeterminate with the philosophy that 
while it would not be problematic to let it be `true', there is no \emph{reason} to
force it to be `true'. More particularly, there is no reason to force it verifiable. This is because it is
not well-founded, so we cannot `dig' the truth value for the formula from any atomic level.
Such an atomic level could be considered unproblematic, like the atomic formulas in
games for first-order logic, but for $C_i\, C_i$ there is no such atomic ground level.
We need to keep on digging, or at least we can never reach a bottom, no matter what.
Also, if the interpretation of $x$ is deleted, then $P(x)$ becomes naturally indeterminate.
Even $x=x$ becomes so with the reading for 
atoms that \emph{it can be directly observed that $x=x$}. Such a reading works 
well also for other first-order atoms.

However, $\mathcal{L}$ is a natural logic for considering partial functions. Indeed, when we
add an isolated element to a model, then all functions become partial. Also, if we delete the
interpretation of a constant symbol, a similar situation is realized. However, 
the game-theoretic semantics nicely facilitates the treatment of such phenomena.
If the game ends up in an atom $R(t_1,\dots , t_n)$, then the following happens.

\begin{enumerate}
\item
The verifier wins if $t_1^{\mathfrak{M},g},\dots , t_n^{\mathfrak{M},g}$ are 
all defined and $R^{\mathfrak{M}}$ is defined on $(t_1^{\mathfrak{M},g},\dots , t_n^{\mathfrak{M},g})$
and we have $(t_1^{\mathfrak{M},g},\dots , t_n^{\mathfrak{M},g}) \in R^{\mathfrak{M}}$.
\item
The falsifier wins if $t_1^{\mathfrak{M},g},\dots , t_n^{\mathfrak{M},g}$ are 
all defined and $R^{\mathfrak{M}}$ is defined on $(t_1^{\mathfrak{M},g},\dots , t_n^{\mathfrak{M},g})$
and we have $(t_1^{\mathfrak{M},g},\dots , t_n^{\mathfrak{M},g}) \not\in R^{\mathfrak{M}}$.
\item
Otherwise the game ends and neither player wins the play of the game.
\end{enumerate}
The treatment of equality atoms $t_1 = t_2$ are similar.

\begin{enumerate}
\item
The verifier wins if $t_1^{\mathfrak{M},g}$ and $t_2^{\mathfrak{M},g}$ are 
defined and we have $t_1^{\mathfrak{M},g} = t_2^{\mathfrak{M},g}$.
\item
The falsifier wins if $t_1^{\mathfrak{M},g}$ and $t_2^{\mathfrak{M},g}$ are 
defined and we have $t_1^{\mathfrak{M},g} \not= t_2^{\mathfrak{M},g}$.
\item
Otherwise the game ends and neither player wins the play of the game.
\end{enumerate}

Generally, we allow terms $t$ in all atoms, while the first take on $\mathcal{L}$ allowed 
first-order variables only. Note that $D$ makes a tuple $(u_1,\dots , u_n)$ not belong to the modified
relation, whether it was undefined or positively true that $(u_1,\dots , u_n)$ is in the 
relation originally. Similarly $I$ positively adds the tuple, no matter the previous status of the tuple.

Finally, the semantics for first-order logic with partially defined entities given in the 
beginning of this document is equivalent to what we obtain for that logic from
the described game-theoretic semantics.

\medskip

\medskip

\medskip

\begin{theorem}
The game-theoretic semantics extends the compositional one on
the first-order fragment.
\end{theorem}

\medskip

\medskip

\medskip

Concerning $\mathcal{L}$ it is of course possible and natural to extend $\mathcal{L}$ so
that it can also modify partial functions $f$ and partially constants $c$ 
the same way it is allowed to modify partial relation symbols. Furthermore, there are
many ways to add $\sim$ to $\mathcal{L}$. Some of the natural choices lead to a 
logic that captures the arithmetic hierarchy. Also, we could, e.g., consider variants of $D$
that only make it undefined whether a tuple belongs to a relation. Indeed, we could naturally
define similar operators for all
other modificationa between positive, undefined and negative instances.

\section{Operators}

We now define a notion of a \emph{modifier} (cf. \cite{gamesandcomp}). This notion works
well with the three-valued logic of partial definitions.

Let $V_R$ denote the full relational vocabulary consisting of 
countably infinitely many relation symbols $R$ for 
each arity $k\in \mathbb{N}$. Similarly, let $V_f$ denote 
the full functional vocabulary consisting of countably 
infintely many function symbols of each arity $k\in \mathbb{Z}_+$.
Moreover, let $V_c$ be the full constant vocabulary consisting of
countably infinitely many constant symbols. Recall
that $\mathrm{VAR}$ denotes a countably infinite set of 
variable symbols. We note that these sets can also be 
considered having higher infinite numbers of symbols, but 
for most purposes, countably infinite sets suffice. 
Let $\mathrm{SYMB}$ denote the 
union \[V_R\cup V_f \cup V_c \cup \mathrm{VAR}.\]

A \textbf{symbol list} is a $4$-tuple of
the form 
\[((R_1,\dots , R_{k_R}), (f_1,\dots , f_{k_f}), 
(c_1,\dots , c_{k_c}), 
(y_1,\dots , y_{k_y}))
\]
that specifies lists of relation symbols,
function symbols, constant symbols and variables. 
Here $k_R,k_f,k_c,k_y\in \mathbb{N}$ are natural numbers. Generally,
these can be any finite natural numbers, including zero.
Nullary relation symbols are allowed. The \textbf{type} of
the above symbol list is the $4$-tuple $(\overline{s}_1,\overline{s}_2,
k_c,k_y)$ such that the following conditions hold. 
\begin{enumerate}
\item
$\overline{s}_1$ is a tuple in $\mathbb{N}^{k_R}$
that gives the arities of the symbols in $(R_1,\dots , R_{k_R})$.
\item
Similarly, $\overline{s}_2$ is a tuple in $\mathbb{Z}_+^{k_f}$
that gives the arities of the symbols in $(f_1,\dots , f_{k_f})$.
\end{enumerate}
%
%
%
%
%
%
%

By a \textbf{symbol set}, we mean a set that 
contains relation symbols, function symbols, constant symbols 
and variables. Thus symbol sets are unions of vocabularies and
variable sets. By an \textbf{interpretation} we mean a
pair $(\mathfrak{M},f)$ where $\mathfrak{M}$ is a model and $f$ is a
variable assignment mapping some set $V\subseteq \mathrm{VAR}$ 
into the domain $M$ of $\mathfrak{M}$. If $\tau$ is a symbol set, we
let $\mathrm{Int}_{\tau}$ denote the class of interpretations $(\mathfrak{M},f)$
where $\mathfrak{M}$ is a model over the vocabulary $\tau \setminus \mathrm{VAR}$
and the domain of $f$ is precisely $\tau \cap \mathrm{VAR}$. These can be
called $\tau$-interpretations.

Let $\sigma$ be a symbol set. A $\sigma$-\textbf{property triple} is a 
product of the form
\[ \mathcal{P}(\mathrm{Int}_{\sigma}) 
\times \mathcal{P}(\mathrm{Int}_{\sigma}) 
\times \mathcal{P}(\mathrm{Int}_{\sigma}) \]
where $\mathcal{P}$ denotes the power set (or power class) operator.
An \textbf{$\ell$-classification} of 
type $(\sigma_1,\dots , \sigma_{\ell})$ is a sequence 
\[(\mathit{Pt}(\sigma_1),  \dots , \mathit{Pt}(\sigma_{\ell}))\]
where each $\mathit{Pt}(\sigma_i)$ is a $\sigma_i$-property triple. 
We let $\mathcal{C}_{(\sigma_1,\dots , \sigma_{\ell})}$ denote the 
class of $\ell$-classifications of type $(\sigma_1 , \dots , \sigma_{\ell})$.

Consider symbol sets $\tau,\sigma_1,\dots , \sigma_{\ell}$, and
let $\overline{\sigma}$ denote $(\sigma_1,\dots , \sigma_{\ell})$, so we
can also write $\mathcal{C}_{\overline{\sigma}}$ for $\mathcal{C}_{(\sigma_1,\dots , \sigma_{\ell})}$. 
Now, a $(\tau,\overline{\sigma})$-\textbf{operator} is a mapping 
\[F_{\tau,
\overline{\sigma}}:\mathrm{Int}_{\tau}
\rightarrow \mathcal{P}(\mathcal{C}_{\overline{\sigma}}) \]
that gives a class of $\ell$-classifications of type $\overline{\sigma}$
for each input $\mathfrak{I}\in \mathrm{Int}_{\mathcal{\tau}}$ 
(recall indeed that $\mathcal{P}$ is the power set---or power class---operator). 
Furthermore, the following natural invariance conditions hold.
\begin{enumerate}
\item
Suppose $\mathfrak{I}\in \mathrm{Int}_{\tau}$
and $\mathfrak{I}'\in \mathrm{Int}_{\tau}$ are isomorphic via an isomorphism $f$. 
Then there is a one-to-one mapping $g$ 
from $F(\mathfrak{I})$ to $F(\mathfrak{I}')$, where we 
let $F$ denote $F_{\tau,\overline{\sigma}}$ for simplicity. 
\item
For each tuple $T = (T_1,\dots , T_{\ell})\in F(\mathfrak{I})$ and
the corresponding tuple $g(T) = (T_1',\dots , T_{\ell}') \in F(\mathfrak{I}')$, 
and for each $i\in \{1, \dots , \ell \}$, there exist 
one-to-one correspondences $h_1: T_i[1] \rightarrow T_i'[1]$,
$h_2: T_i[2] \rightarrow T_i'[2]$ and $h_3: T_i[3] \rightarrow T_i'[3]$,
where $X[j]$ denotes the $j$th member of the triple $X$. 
\item 
For each of these one-to-one correspondences $h: T_i[j] \rightarrow T_i'[j]$ 
and for each $\mathfrak{I}''$ and $\mathfrak{I}'''$ such
that $h(\mathfrak{I}'') = \mathfrak{I}'''$, there exists an
isomorphism $r:\mathfrak{I}''\rightarrow \mathfrak{I}'''$. 
\end{enumerate}

The above conditions simply force isomorphism invariance.

Now, let $f:\mathrm{SYMB} \rightarrow \mathrm{SYMB}$ be a
bijection such that for each $V\in \{V_R, V_f, V_c, \mathrm{VAR}\}$,
the restriction $f\upharpoonright V$ is a
bijection from $V$ onto $V$ that preserves arities of
relation symbols and function symbols. We call such a bijection a
\textbf{symbol name permutation}. 
Given a symbol name permutation $f$, by $f(\tau)$ we denote the
symbol set $\{f(x)\ |\ x\in \tau \}$ obtained 
from the symbol set $\tau$ in the natural way. 
For $\overline{\sigma} = (\sigma_1,\dots , \sigma_{\ell})$, 
we let $f(\overline{\sigma})$ 
denote $(f(\sigma_1), \dots , f(\sigma_{\ell}))$. 
An \textbf{$f$-isomorphism}
from $(\mathfrak{M},g)\in \mathrm{Int}_{\tau}$ to 
$(\mathfrak{N},h)\in \mathrm{Int}_{f(\tau)}$ is a
bijection $p$ from the domain $M$ of $\mathfrak{M}$ 
onto the domain $N$ of $\mathfrak{N}$ such that the
following conditions hold.
\begin{enumerate}
\item
The domains of the assignments $g$ and $h$ have 
the same number of variables, and for each $x$ in
the domain of $g$, we have $g(x) = m\in M$ if and 
only if $h(f(x)) = p(m)$.  
\item
The vocabularies of $\mathfrak{M}$ and $\mathfrak{N}$ are of
the same type in the sense that the following conditions hold.
\begin{enumerate}
\item
For each $k\in\mathbb{N}$, both vocabularies have the same number of
relation symbols of arity $k$.
\item
For each $k\in\mathbb{Z}_+$, both vocabularies have the same number of
function symboss of arity $k$.
\item
Both vocabularies have the same number of constant symbols. 
\end{enumerate}
\item
The following is true for each constant symbol $c\in \tau$
and the corresponding symbol $d := f(c)$. It holds that $c^{\mathfrak{M}}$ is
defined to be $m\in M$ if and only if $d^{\mathfrak{N}}$ is
defined to be $p(m) \in N$. 
\item
For each arity $k$, the following is true for each relation symbol $R$ of arity $k$
and the corresponding symbol $S := f(R)$ and for each tuple 
$(m_1,\dots , m_k)\in M^k$ of domain elements. 
\begin{enumerate}
\item
$R^{\mathfrak{M}}$ is defined on the tuple $(m_1,\dots , m_k)$ if
and only if $S^{\mathfrak{N}}$ is
defined on the tuple $(p(m_1), \dots , p(m_k))$.
\item
$(m_1,\dots , m_k)$ is defined to be in $R^{\mathfrak{M}}$ if
and only if $(p(m_1),\dots , p(m_k))$ is
defined to be in $S^{\mathfrak{N}}$. Note that this, together
with the above condition, implies that $(m_1,\dots , m_k)$ is
defined to not be in $R^{\mathfrak{M}}$ if
and only if $(p(m_1),\dots , p(m_k))$ is
defined to not be in $S^{\mathfrak{N}}$. 
\end{enumerate}
\item
For each arity $k$, the following holds for each relation symbol $q$ of arity $k$
and the corresponding symbol $r := f(q)$ and for each tuple 
$(m_1,\dots , m_k)\in M^k$. 
\begin{enumerate}
\item
$q^{\mathfrak{M}}$ is defined on the tuple $(m_1,\dots , m_k)$ if
and only if $r^{\mathfrak{N}}$ is
defined on the tuple $(p(m_1), \dots , p(m_k))$.
\item
$q^{\mathfrak{M}}(m_1,\dots , m_k)$ is defined to be $m\in M$ if
and only if $r^{\mathfrak{N}}(p(m_1),\dots , p(m_k))$ is
defined to be $p(m)$. We note that condition $a$ actually follows 
already from condition $b$. 
\end{enumerate}
\end{enumerate}

Let $f$ be a symbol name permutation. 
Now, an $f$-\textbf{renaming} of a $(\tau,\overline{\sigma})$-operator is an $(f(\tau),
f(\overline{\sigma}))$-operator $F$ such that the following condition holds.

\begin{enumerate}
\item
Suppose there exits an $f$-isomorphims from $\mathfrak{I}\in
\mathrm{Int}_{\tau}$ to $\mathfrak{I}' \in \mathrm{Int}_{f(\tau)}$. 
Then there is a one-to-one mapping $g$ 
from $F(\mathfrak{I})$ to $F(\mathfrak{I}')$.
\item
For each tuple $T = (T_1,\dots , T_{\ell})\in F(\mathfrak{I})$ and
the corresponding tuple $g(T) = (T_1',\dots , T_{\ell}') \in F(\mathfrak{I}')$, 
and for each $i\in \{1, \dots , \ell \}$, there exist 
one-to-one correspondences $h_1: T_i[1] \rightarrow T_i'[1]$,
$h_2: T_i[2] \rightarrow T_i'[2]$ and $h_3: T_i[3] \rightarrow T_i'[3]$,
where $X[j]$ denotes the $j$th member of the triple $X$.
\item 
For each of these one-to-one correspondences $h: T_i \rightarrow T_i'$ 
fixed above and for each $\mathfrak{I}''$ and $\mathfrak{I}'''$ such
that $h(\mathfrak{I}'') = \mathfrak{I}'''$, there exists an $f$-isomorphism $r:\mathfrak{I}''\rightarrow \mathfrak{I}'''$. 
\end{enumerate}

Let $\overline{S} = (\overline{S}_{1}, \dots , \overline{S}_{\ell})$ be a
tuple of $\ell$ symbol list types, so each $\overline{S}_i$ is a 
symbol list type. A \textbf{symbol list tuple} $\overline{\rho} = (\overline{\rho}_1,\dots ,
\overline{\rho}_{\ell})$ is a tuple of symbol lists $\overline{\rho}_i$.
The symbol list tuple $\overline{\rho} = (\overline{\rho}_1,\dots ,
\overline{\rho}_{\ell})$ is of type $\overline{S} = (\overline{S}_{1}, \dots , \overline{S}_{\ell})$ if
each $\overline{\rho}_i$ has the type $\overline{S}_i$. 
Now, an \textbf{operator} $\dot{O}$ of type $\overline{S}$ is a mapping that
takes as input a 
symbol set $\tau$ and a symbol list tuple $\overline{\rho} = (\overline{\rho}_{1},
\dots , \overline{\rho}_{\ell})$ of type $\overline{S}$, and based on this information, outputs a $(\tau,\overline{\sigma})$-operator $\dot{O}_{\tau}(\overline{\rho})$ 
such that the following conditions hold. 
\begin{enumerate}
\item
Note first indeed 
that the $(\tau, \overline{\sigma})$-operator is denoted by $\dot{O}_{\tau}(\overline{\rho})$ where we have
the symbol list tuple $\overline{\rho}$ instead of the symbol set tuple $\overline{\sigma}$. The symbol set tuple $\overline{\sigma}$ has $\ell$
symbol sets, so we can 
write $\overline{\sigma} = (\sigma_1,\dots , \sigma_{\ell})$.  
\item
For each ${\sigma_i}$ in  
the tuple $\overline{\sigma} = (\sigma_1,\dots , \sigma_{\ell})$, we 
have $\sigma_i \subseteq \tau \cup \mathit{symb}(\overline{\rho}_i)$
where $\mathit{symb}(\overline{\rho}_i)$ denotes the set of
all symbols that appear in the tuples of 
the symbol list $\overline{\rho}_i$. 
\end{enumerate}
%
%
%
%
%
%
%
%
%
%
%
%

%
%
%
Note that $\tau$ in the input of an operator can vary freely, 
and so can $\overline{\rho}$, as long as it is of type $\overline{S}$. 
An operator of \textbf{relaxed type} $\overline{S}$ lifts the 
restrictions $\sigma_i\subseteq \tau \cup \mathit{symb}(\overline{\rho}_i)$, i.e.,
those restrictions do not apply. Nevertheless, we give $\overline{\rho}$ of
type $\overline{S}$ (together
with $\tau$) as a
convenient syntactic input to an operator of
relaxed type $\overline{S}$. An operator of \textbf{$\ell$-relaxed type}
allows any symbol list tuple $\overline{\rho}$ with $\ell$
symbol lists as an 
input, so $\overline{\rho}$ does not have to be of type $\overline{S}$, 
and again the restrictions $\sigma_i
\subseteq \tau \cup \mathit{symb}(\overline{\rho}_i)$ do not apply.

Now, let $\tau$ and $\tau'$ be a symbol sets.
Let $\overline{\rho} = (\overline{\rho}_1, \dots , \overline{\rho}_{\ell})$
and $\overline{\rho}' = (\overline{\rho}_1',\dots , \overline{\rho}_{\ell}')$ be
symbol list tuples. We say
that $(\tau , \overline{\rho})$ and $(\tau' , \overline{\rho}')$ are
\textbf{naming variants} if there is a symbol name permutation $f$ 
such that $(f(\tau) , f(\overline{\rho})) = (\tau' , \overline{\rho}')$, 
where $f(\overline{\rho})$ denotes the symbol list tuple obtained
from $\overline{\rho}$ by replacing every
symbol $\# \in \mathrm{SYMB}$ that appears in $\overline{\rho}$ by the
corresponding symbol $f(\#)$. 
The restriction of $f$ to the set of symbols 
appearing in $(\tau , \overline{\rho})$ is 
then the corresponding 
\textbf{renaming bijection} from $(\tau , 
\overline{\rho})$ to $(\tau' , \overline{\rho}')$.

An operator $\dot{O}$ is \textbf{renaming invariant} if
for any inputs $(\tau,\overline{\rho})$ and $(\tau', \overline{\rho}')$ 
that are naming variants with a renaming bijection $b$ 
from $(\tau , 
\overline{\rho})$ to $(\tau' , \overline{\rho}')$, 
the output operator $\dot{O}_{\tau'}(\overline{\rho}')$ is an $f$-renaming of 
$\dot{O}_{\tau}(\overline{\rho})$ for some $f$ that extends $b$. Renaming 
invariant operators are logically quite natural, whether they are of
some type $\overline{S}$, or of relaxed type $\overline{S}$ or 
just of $\ell$-relaxed type. 

Now, consider an operator $\dot{O}$ of type $\overline{S} = (\overline{S}_1, 
\dots , \overline{S}_{\ell})$, and let $O$ be a 
related symbol referring to $\dot{O}$. We continue treating operators of
relaxed type $\overline{S}$ and $\ell$-relaxed type
similarly, so alternatively $\dot{O}$ can be an 
operator of one of the relaxed types. 
We define a semantics for formulae $O\overline{\rho}(\varphi_1,\dots , \varphi_{\ell})$, where $\overline{\rho}$ is a 
symbol list tuple of type $\overline{S}$ (or of 
any type with $\ell$ symbol lists in the case of an operator of
$\ell$-relaxed type). Let $(\mathfrak{M},g)$ be a $\tau$-interpretation.
We define that $\mathfrak{M},g\models^+ O\overline{\rho}(\varphi_1,
\dots , \varphi_{\ell})$ if and
only if there exists a tuple $C$ in $F((\mathfrak{M},g))$,
where $F$ denotes $\dot{O}_{\tau}(\overline{\rho})$, such
that the following conditions hold. 
\begin{enumerate}
\item 
For each $i\in \{1,\dots , \ell\}$, let $C_i$ denote
the $i$th triple in $C$, and for each $j\in\{1,2,3\}$,
let $C_i[j]$ denote the $j$th member of $C_i$. 
Then, for each $i\in \{1, \dots , \ell \}$, we have  
$\mathfrak{N},h\models^+ \varphi_i$ for
each $(\mathfrak{N},h)\in C_i[1]$. 
\item 
Similarly, for each $i\in \{1, \dots , \ell \}$, we
have $\mathfrak{N},h \models^- \varphi_i$
for each $(\mathfrak{N},h)\in C_i[3]$.  
\item
For each $i\in \{1, \dots , \ell \}$, we
have $\mathfrak{N},h\not\models^+ \varphi_i$
and $\mathfrak{N},h\not\models^-\varphi_i$
for each $(\mathfrak{N},h)\in C_i[2]$.
\end{enumerate}

Note that renaming invariant modifiers are insensitive to 
change of names of operators. This can be natural. Further 
invariance conditions are also natural for various purposes. It is
easy to generalize, e.g., to model sets. There the input is a
set (or class) of interpretations $(\mathfrak{M},g)$. 
However, also simpler operators are natural. A \textbf{simple} 
operator is a renaming invariant
operator of type $\overline{S} = ((\emptyset,\emptyset,
\emptyset,\emptyset))$ (where $\emptyset$ is
the empty tuple)
which returns, given $(\mathfrak{M},g)$, an 
interpretation $(\mathfrak{N},h)$ that has the same
symbol set as $(\mathfrak{M},g)$. A \textbf{very simple} 
operator is a simple operator that returns,
given $(\mathfrak{M},g)$, an interpretation $(\mathfrak{N},g)$ 
where the assignment has not changed and the
domain of $\mathfrak{N}$ is the same as that of $\mathfrak{M}$. 
Simple and very simple operators come with extensions that 
deal with several formulae instead of only one.

Note that the
considered three-valued logic is of course just a generalization of 
the two-valued approach, and we can always define operators
where each of the second sets $C_i[2]$ is always empty. That
covers the two-valued approach. Finally, (for example) the 
standard connectives and generalized quantifiers are
easy to deal with via modifiers, but it is perhaps more interesting to
consider `more dynamic' model modifiers not
relating directly to quantifier statements or plain connectives.

\end{document}